\documentclass{amsart}
\usepackage{enumerate}
\usepackage{amssymb,amsmath}
\usepackage{latexsym}
\usepackage{graphicx}
\usepackage[]{color}
\usepackage{ifthen}

\newcommand{\cA}{\mathcal{A}}
\newcommand{\cQ}{\mathcal{Q}}
\newcommand{\R}{\mathbb{R}}

\newtheorem{theorem}{Theorem}[section]
\newtheorem{lemma}[theorem]{Lemma}
\newtheorem{proposition}[theorem]{Proposition}
\newtheorem{corollary}[theorem]{Corollary}
\newtheorem{problem}{Problem}[section]
\newtheorem{conjecture}{Conjecture}[]
\newtheorem{question}{Question}[section]

\theoremstyle{definition}
\newtheorem{definition}[theorem]{Definition}

\newtheorem{examples}[theorem]{Examples}

\theoremstyle{remark}

\numberwithin{equation}{section}
\newif\ifmargin
\margintrue
\newif\ifskip
\skiptrue

\newcommand{\NN}{\mathbb{N}}
\newcommand{\PP}{\mathbb{P}}

\newcommand{\EE}{\mathbb{E}}

\newcommand{\cP}{\mathcal{P}}
\newcommand{\cD}{\mathcal{D}}

\newcommand{\cH}{\mathcal{H}}
\newcommand{\bP}{\mathbf{P}}

\newcommand{\goesto}{\rightarrow}

\begin{document}

\title{Almost unimodal  and real-rooted \\ graph polynomials}

\author{Johann A. Makowsky}
\address{Department of Computer Science, Technion - IIT, Haifa, Israel}
\email{janos@cs.technion.ac.il}

\author{Vsevolod Rakita} 
\address{Department of Mathematics, Technion - IIT, Haifa, Israel}
\email{vsevolod@campus.technion.ac.il}



\date{\today}


\begin{abstract}
It is well known that the coefficients of the matching polynomial are unimodal.
Unimodality of the coefficients (or their absolute values) of other graph polynomials
has been studied as well.
One way to prove unimodality is to prove real-rootedness.

Recently I. Beaton and J. Brown (2020) proved the for almost all graphs the coefficients of
the domination polynomial form a  unimodal sequence, and
C. Barton, J. Brown and D. Pike (2020) proved that the forest polynomial (aka acyclic polynomial)
is real-rooted iff $G$ is a forest.

Let $\cA$ be a graph property, and let $a_i(G)$ be the number of induced subgraphs of order $i$ of a graph $G$
which are in $\cA$.
Inspired by their results we prove:

{\bf Theorem:} If $\cA$ is the complement of a hereditary property, then for almost all graphs in $G(n,p)$
the sequence $a_i(G)$ is unimodal.

{\bf Theorem:} If $\cA$ is a hereditary property which contains a graph which is not a clique or the complement of a clique,
then the graph polynomial $P_{\cA}(G;x) = \sum_i a_i(G) x^i$ is real-rooted iff $G \in \cA$.
\end{abstract}

\maketitle
\tableofcontents

\section{Introduction}
\subsection{Graph polynomials and their properties}
A graph polynomial is a graph invariant $P(G)$  with values in a polynomial ring, 
usually a subring of $\R[\bar{x}]$. $P(G)$ is univariate if for every graph $G$ the
polynomial has one indeterminate. Particular graph polynomials, such as the chromatic polynomial,
the characteristic polynomial,
the matching polynomial, and the Tutte polynomial have been studied extensively.
In \cite{makowsky2008zoo} the first author initiated a project to develop a comparative study of graph polynomials.
Its purpose is to find properties of graph polynomials which are shared by infinitely many uniformly described
families of graph polynomials.
A graph property $\cA$ is called hereditary if it is closed under taking induced subgraphs, i.e. if $G\in \cA$ and $H$ is an induced subgraph of $G$, then $H\in \cA$. A graph property $\cA$ is called co-hereditary if it is the complement of a hereditary graph property.
In this paper we present two general theorems for univariate graph polynomials $P_{\cA}(G;x)$ 
which are generating functions
of hereditary, respectively co-hereditary graph properties $\cA$. We show that for  $\cA$
co-hereditary, the graph polynomial $P_{\cA}(G;x)$ is unimodal for almost all graphs,
and that for $\cA$ hereditary which contains at least one graph which is not a clique or a complement of a clique,
$P_{\cA}(G;x)$ is real-rooted iff $G \in \cA$.

\subsection{Basic Definitions and Notation}
Throughout this paper, except where stated otherwise, we consider undirected simple graphs, i.e. graphs without parallel edges or loops, with labelled vertices. For a graph $G$, we denote by $V(G)$ its vertex set and by $E(G)$ its edge set. We denote by $n(G)$ the order of $G$, i.e $n(G)=|V(G)|$, and by $e(G)$ the size of $G$, i.e. $e(G)=|E(G)|$.
Let $G$ be a graph, and let $U\subseteq V(G)$ be a subset of vertices of $G$. The graph {\em induced} by $U$ in $G$, denoted $G[U]$, is the graph with vertex set $U$ where two vertices age incident if and only if they are incident in $G$.

For graphs $G,H$, we say that $H$ is a {\em subgraph} of $G$ if $H$ can be obtained from $G$ by deleting vertices and edges, and we say that $H$ is an {\em induced subgraph} of $G$ if there is a set $U\subseteq V(G)$ such that $H$ is isomorphic to $G[U]$. If $H$ is an induced subgraph of $G$, we write $H<G$.
A {\em graph property}  $\cA$ is a family of graphs closed under isomorphism, i.e. if for a graph $G$, $G\in \cA$ and $H$ is a graph isomorphic to $G$, then $H\in \cA$. 
A graph property $\cA$ is said to be {\em non-trivial} if it is not empty, and there is a graph $G$ such that $G \not \in \cA$.
We will consider certain types of graph properties:
\begin{definition}
A graph property $\cA$ is called {\em hereditary} if it is closed under taking induced subgraphs, i.e. if $G\in \cA$ and $H$ is an induced subgraph of $G$, then $H\in \cA$. A graph property $\cA$ is called {\em co-hereditary} if it is the compliment of a hereditary property, i.e. if there exists a hereditary property $\cA'$ such that for all graphs $G$, $G\in \cA$ if and only if $G\not \in \cA$.
\end{definition}

For $n\in \NN$ and $p\in[0,1]$, a {\em random graph} $G\in \mathcal{G}(n,p)$ is a graph with $n$ vertices, where every two vertices are incident with probability $p$, independently of the others.

We say that a statement holds for {\em almost all graphs} if the proportion of isomorphism classes of graphs for which the statement holds tends to 1 as the order of the graphs tends to infinity, that is if we denote by $U(n)$ the collection of all isomorphism classes of graphs of order $n$ for which the statement holds, and by $\mathcal{G}(n)$ the collection of all isomorphism classes of graphs of order $n$, then $\lim_{n \goesto \infty} \frac{U(n)}{\mathcal{G}(n)}=1$. Equivalently, a statement holds for almost all graphs if the probability that it holds for a random graph $G\in \mathcal{G}(n,1/2)$ tends to 1 as $n$ tends to infinity.

\subsection{Real-rooted and unimodal graph polynomials}

Let $F(x) \in \R[x]$ be an univariate polynomial of degree $d$ with real coefficients, 
$$F(x) = \sum_{i=0}^d a_i x^i$$
\begin{enumerate}[(i)]
\item
$F(x)$ is {\em real-rooted} if all its roots are in $\R$.
\item
The coefficients of
$F(x)$ are {\em log-concave} if for all $1 \leq j \leq d-1$
\\
{\em $a_j^2 \geq a_{j-1} a_{j+1}$.} 
\item
The coefficients of
$F(x)$ are {\em unimodal with mode $k$} if 
\\
{\em $a_i \leq a_j$ for $ 0\leq i < j \leq k$}
and
{\em $a_i \geq a_j$ for $ k \leq  i < j \leq d$.}
\item
$F(x)$ is {\em absolute unimodal with mode $k$ (log-concave)} if 
\\
the absolute values of $a_i$
are unimodal (log-concave).  
\item
These definitions, except for (i), apply to any sequence $a_i, 0 \leq i \leq d$,
even if it not interpreted as a sequence of coefficients of a polynomial.
\end{enumerate}

\begin{theorem}[Folklore]
\label{thm:NewtonsTheorem}
(i) implies (ii), (ii) implies (iii) and
none of the reverse implications holds.
\end{theorem}
The first part of the theorem is is implied by Newton's theorem. For a proof one may consult \cite{branden2015unimodality}.

The sequence ${n \choose i}$ of the number of subsets of order $i$ of $V(G)$ is log-concave, hence unimodal.
More interestingly, 
let $m_i(G)$
the sequence of the number of edge independent subsets (matchings) of $E(G)$ of order $i$.
The numbers $m_i(G)$ are also the coefficient of the generating matching polynomial
$$
M(G;x) =\sum_i m_i(G) x^i.
$$
\begin{theorem}
$M(G;x)$ is real-rooted, hence unimodal. 
\end{theorem}
There are two independent proofs of this.
It follows from the fact that all the roots of $M(G;x)$ are real for all graphs $G$, \cite{heilmann1970monomers},
see also \cite{gutman2016survey},
using Theorem \ref{thm:NewtonsTheorem}. Unimodality was also shown directly by A. Schwenk, \cite{schwenk1981unimodal}. Additionally, C

Let $in_i(G)$, $0\leq i \leq n(G)$,
the sequence of the number of vertex independent subsets of $V(G)$ of order $i$. Denote by
$I(G;x) = \sum_i in_i(G)x^i$
the independence polynomial of $G$.
Real-rootedness, and unimodality of $I(G;x)$ has been studied extensively, see e.g. \cite{zhu2020unimodality}\cite{brown2018unimodality}\cite{cutler2017maximal} for some recent results and \cite{levit2005independence}\cite{levit2006independence} for a general introduction.

\begin{theorem}
\begin{enumerate}[(i)]
\item
$I(G;x)$ is not unimodal, 
\cite{alavi1987vertex}.
However, it is easily seen that the set of counterexamples $G$
given in that paper
has measure $0$ among the random graphs $\mathcal{G}(n,p)$.
\item
For claw-free graphs the sequence of coefficients of $I(G;x)$ is real-rooted, hence unimodal, 
\cite{chudnovsky2007roots,bencs2014christoffel}.
It is easily seen that the claw-free graphs have measure $0$ among the random graphs $G(n,p)$.
\end{enumerate}
\end{theorem}

This leaves open whether $I(G;x)$ is real-rooted, or at least unimodal, for other graph classes. 
Specifically, we may ask whether $I(G;x)$ is unimodal for ``most graphs'', in the following sense:

\begin{definition}
Let $P$ be a graph polynomial. We say $P$ is {\em almost unimodal} if for almost all graphs $G$ the polynomial 
$P(G;x)$ is unimodal. In other words, $P$ is almost unimodal if 
for random graphs $G \in \mathcal{G}(n, 1/2)$ we have
$$
\lim_{n\goesto \infty} \PP(P(G(n,1/2),x) \text{ is unimodal} )=1.
$$
\end{definition}

\begin{problem}
Is $I(G;x)$  almost unimodal?
\end{problem}

Let $\chi(G,x) = \sum_i c_i(G) x^i$ be the chromatic polynomial of $G$.
The case of the chromatic polynomial of a graph is slightly different.
The sequence $c_i(G)$ is alternatingly positive and negative. However,
it was conjectured by R.C. Read, \cite{read1968introduction}, that the absolute values 
$|c_i(G)|$ form a unimodal sequence.
J. Huh, \cite{huh2012milnor} finally proved the conjecture.
\begin{theorem}[J. Huh, 2012]
For every graph $G$ the chromatic polynomial $\chi(G,x)$ is absolute unimodal.
In fact the sequence $|c_i(G)|$ is log-concave.
\end{theorem}

\subsection{Counting induced subgraphs of a graph}

Both $m_i(G)$ and $in_i(G)$ are graph parameters counting subgraphs of $G$ which satisfy a graph property $\cA$,
graphs of degree exactly $1$, and edge-less graphs respectively.

In this paper,
we are interested in the question for which properties $\cA$ the corresponding counting parameters are real-rooted, 
unimodal, or almost unimodal.

Given  a graph property $\cA$,
and denote by
\begin{gather}
c^{\cA}_i(G) = | \{ A \subseteq V(G): |A|=i, G[A] \in \cA \}| \notag\\
P_{\cA}(G;x) =
\sum_i c^{\cA}_i(G) x^i =
\sum_{S \subseteq V(G),G[S]\in \cA} x^{|S|} \notag
\end{gather}
the number of induced subgraphs of $G$ in $\cA$ of order $i$ and its corresponding graph polynomial.
For each $G$ the polynomial
$P_{\cA}(G;x)$ is the generating function of $\cA$.

We now consider a slight generalization of graph properties:
\begin{definition}
Let $\cP$ be a family of pairs $(G,S)$ where $G$ is a graph and $S\subseteq V(G)$ is a set of vertices. We say that $\cP$ is an {\em augmented graph property} if it is closed under $S$ preserving graph isomorphisms, i.e. if $(G,S)\in \cP$, and $f:G \goesto H$ is a graph isomorphism, then $(H,f(S))\in \cP$ (where $f(S)=\{f(v):v\in S\}$).
\end{definition}

We can define graph parameters satisfying an augmented graph property, similar to regular graph parameters:
\begin{definition}
Let $\cP$ be an augmented graph property, and let $G$ be a graph. Denote by
$$
c^{\cP}_i(G) = | \{ S \subseteq V(G): |S|=i, (G, S) \in \cP \}|
$$
the number of subsets $S \subseteq V(G)$ 
of size $i$ with $(G, S) \in \cP$, and by
$$
F_{\cP}(G;x) =
\sum_i c^{\cP}_i(G) x^i =
\sum_{S \subseteq V(G), (G, S) \in \cP} x^{|S|}
$$
its corresponding graph polynomial.

\end{definition}
We will be interested in a certain type of augmented graph property:
\begin{definition}
Let $\cP$ be an augmented graph property.
We say $\cP$ is {\em upward monotone}, 
if whenever $(G, S) \in \cP$ and $S\subseteq S' \subseteq V(G)$, 
then $(G, S') \in \cP$.
\end{definition}

Typical examples which we will use in this paper are:
\begin{examples}
\phantom{hjg}
\begin{enumerate}[(i)]
\item
Recall that a {\em dominating set} in a graph $G$ is a set $D\subseteq V(G)$ such that every vertex $v\in V(G)$ is either in $D$ or is incident to a vertex in $D$.
The augmented graph property $DOM$ consists of all graphs $(G, D)$ with a distinguished dominating set $D \subseteq V(G)$.
This is an upward monotone augmented graph property.
The corresponding graph polynomial is the {\em domination polynomial}, $F_{DOM}(G;x) = DOM(G;x)$.

\item Let $G$ be a graph. A {\em zero forcing} set in $G$ is a set $S\subseteq V(G)$ such that there is an ordering $\{v_1,v_2,...,v_k\}$ of the vertices in $V(G)-S$ with the property that for each $i$, there is a vertex in the set $S\cup \{v_1,...,v_{i-1}\}$ that is incident to $v_i$ and is not incident to $v_j$ for all $j>i$. The augmented graph property $ZF$ consists of all graphs $(G,S)$ together with a distinguished zero forcing set $S$. This is an upward monotone augmented graph property.
The corresponding graph polynomial is the {\em zero forcing } polynomial. See \cite{boyer2019zero}.

\item
Given a graph property $\cA$, we can define an associated augmented graph property $\cP_{\cA}$
by 
$$
\cP_{\cA} = \{(G, S) :  G[S] \in \cA \}.
$$
For every graph property $\cA$  we have
$$
c^{\cA}_i(G) = c^{\cP_{\cA}}_i(G) \text{   and  } P_{\cA}(G;x)= F_{\cP_{\cA}}(G;x)
$$
The converse is not true, as shown in Theorem \ref{th:xx} below, with a proof given at the end of this section.
\end{enumerate}
\end{examples}

\begin{theorem}[\cite{makowsky2019logician}]
\label{th:MRK}
\label{th:xx}
There is no graph property $\cA$ such that
for all graphs $G$  and for all $i \leq n(G)$ we have
$$
c^{DOM}_i(G) = c^{\cA}_i(G).
$$
\end{theorem}
For convenience of the reader we include here the proof.
\begin{proof}
We first compute some values for $c^{DOM}_i(G)$ the graphs $K_2$ and its complement graph $\bar{K_2} = E_2$. 
\begin{gather}
c^{DOM}_1(K_2) = 2, c^{DOM}_1(E_2) = 0 \label{dom}
\end{gather}
Now, assume, for contradiction, there is such a $\cA$. 
We distinguish two cases.
\\
Case 1: $K_1 \in \cA$.
\\
Then 
$c^{\cA}_1(E_2) = 2$ because $K_1 \in \cA$.
However, 
$c^{DOM}_1(E_2) = 0$
by equation (\ref{dom}), a contradiction.

Case 2: $K_1 \not \in \cA$.
Then
$c^{\cA}_1(K_2) = 0$ because $K_1 \not \in \cA$.
However, 
$c^{DOM}_1(K_2) = 2$
by equation (\ref{dom}), a contradiction.
\end{proof}

\subsection{Unimodality for almost all graphs}

I. Beaton and J. Brown, \cite{beaton2020unimodality}, very recently proved the following theorem.

\begin{theorem}[\cite{beaton2020unimodality}]
Let $d_i(G)$ be the number of dominating sets $D$ of a graph $G$ with $|D| =i$.
For almost all graphs $G$ the sequence $d_i(G), 0\leq i \leq n(G)$ is unimodal. 
\end{theorem}

Their proof suggests the following conjecture:

\begin{conjecture}
\label{conjecture}
Let $\cP$ be an upward monotone augmented graph property.
For a graph $G$, denote by $c^{\cP}_i(G)$ be the number of subsets $S\subseteq V(G)$ with $|S| =i$ such that
$(G, S) \in \cP$.
For almost all graphs $G$ the sequence $c^{\cP}_i(G), 0\leq i \leq n(G)$ is unimodal. 
\end{conjecture}

Our main result here is a first step in proving this conjecture.

\begin{theorem}[Almost Unimodality Theorem]
\label{th:main-1}
Let $\cA$ a non-trivial co-hereditary graph property.
Let $c^{\cA}_i(G)$ be the number of subsets $S\subseteq V(G)$ of a graph $G$ with $|S| =i$ such that
$G[S] \in \cA$.
For almost all graphs $G$, the sequence $c^{\cA}_i(G), 0\leq i \leq n(G)$ is unimodal. 
\end{theorem}

Theorem \ref{th:main-1} does not imply Conjecture \ref{conjecture} because of Theorem \ref{th:xx}.

\subsection{Real-rooted graph polynomials}

In \cite{barton2020acyclic}\footnote{
In the paper \cite{barton2020acyclic}
$F(G;x)$ is called the {\em acyclic polynomial}
This is an unfortunate choice, as the acyclic polynomial exists in the literature
as one of the version the matching polynomial, also called the defect matching polynomial,  $\mu(G;x)$,
used originally in \cite{heilmann1970monomers}.
}.
the following is shown:
\begin{theorem}[\cite{barton2020acyclic}]
$FOR(G;x)$ is real-rooted iff $G$ is a forest.
\end{theorem}

It turns out that their proof (almost verbatim) generalizes
to the following:

\begin{theorem}[Real-rootedness Theorem]
\label{th:main-2}
Let $\cA$ be hereditary and with a graph $G \in \cA$ which is neither a clique nor an edgeless graph.
Then $P_{\cA}(G;x)$ is real-rooted iff $G \in \cA$.
\end{theorem}

For $\cA$ the class of edgeless graphs, $P_{\cA}(G;x)$ is the independence polynomial $I(G;x)$.
Theorem \ref{th:main-2} can not be extended to cover $I(G;x)$, because both the acyclic and the generating matching polynomials
are real-rooted by \cite{heilmann1970monomers}.
Let $g(G;x)$ be the generating matching polynomial of $G$. If $L(G)$ is the line graph of $G$, then
$I(L(G);x) = g(G;x)$. Therefore, $In(G;x)$ is real-rooted for line graphs, 
and by 
\cite{chudnovsky2007roots}
also for claw-free graphs.
The clique polynomial $Cl(G;x)$ is the independence of the complement graph $G^c$,
$Cl(G;x) = I(G^c;x)$, hence the theorem also fails for the Clique polynomial.


\section{Proofs of The Almost Unimodality Theorem \ref{th:main-1}} 

In this section we prove our main results. We begin by proving
Corollary \ref{cor}, which asserts
unimodality of $P_\cA$ for a co-hereditary graph property $\cA$ under the condition that
\begin{gather}
\frac{c^{\cA}_i}{{n \choose i }}\geq \frac{n-i}{i+1}. 
\tag{*}
\label{cond}
\end{gather}
Then we use Theorem \ref{JLR} due to S. Janson, T. Luczak and A. Ruci\'nski, \cite{janson1990exponential,janson2011random}, 
to show that condition (\ref{cond}) is satisfied for almost all graphs.

\subsection{A criterion for unimodality}
Throughout this subsection, we write $c_k$ for $c_k^{\cA}$.
\begin{lemma}
\label{lemma1}
Let $\cA$ be a co-hereditary graph property. 
then for $0\leq k <n/2$, $c_k \leq c_{k+1}$.
\end{lemma}
%

\begin{proof}
Fix $k<n/2$. If $c_k=0$, the claim is trivial, so assume $c_k>0$.
Denote $D_{k+1}=\{B\subseteq V(G):|B|=k+1,\exists A\in C_k, A\subseteq B \}$. Note that if $G[A]\in \cA$, then $G[A\cup\{v\}]\in \cA$ for every vertex $v\in V(G)$, so $D_{k+1}\subseteq C_{k+1}$, and hence it is enough to prove $c_k\leq |D_{k+1}|$.

Consider the bipartite graph $(X\cup Y,E)$ where $X=C_k$, $Y=D_{k+1}$ and there is an edge between $A\in X$ and $B\in Y$ if and only if $A\subseteq B$. Note that every vertex in $X$ has degree $n-k$, so there are $c_k(n-k)$ edges. On the other hand, if $|Y|<|X|$, there is a vertex $B\in Y$ with degree larger then $n-k$, but the degree of a vertex in $Y$ is at most $k+1$, so we have $n-k<k+1$, but then $n/2\leq k$, which is a contradiction.

Overall, we have $c_k=|X|\leq |Y|=|D_{k+1}|$ as required.
\end{proof}
A specialized version of this lemma was proved for dominating sets in \cite{beaton2020unimodality}, and for zero forcing sets in \cite{boyer2019zero}.

\begin{lemma}
\label{lemma2}
Let $\cA$ be as before, $G$ a graph of order $n$, and $k\geq n/2$. If $\frac{c_k}{{n \choose k }}\geq \frac{n-k}{k+1}$, 
then $c_i\geq c_{i+1}$ for $i\geq k$. 
\end{lemma}
\begin{proof}
Denote $r_i=\frac{c_i}{{n \choose i}}$. Note that for all $i$, $r_{i+1}\geq r_i$: if we denote by $A_{i+1}=\{(v,S):v\in S, S\in C_{i+1}\}$ and by $B_i=\{(v,S): v\in V-S, S\in C_i\}$, we have an injective mapping $f:B_i \goesto A_{i+1}$ defined by $f(v,S)=(v,S\cup\{v\})$. Thus, we have that 
$$
(n-i)c_i=|B_i|\leq |A_{i+1}|=(i+1)c_{i+1}
$$
and so
$$
r_{i+1}=\frac{c_{i+1}}{{n \choose i+1}}\geq \frac{n-i}{i+1} \frac{c_i}{{n \choose i+1}}=\frac{c_i}{{n \choose i}}=r_{i}
$$
Now, if $r_k\geq \frac{n-k}{k+1}$, for $i\geq k$ we have $$r_i\geq r_k \geq \frac{n-k}{k+1} \geq \frac{n-i}{i+1}$$ and so
$$\frac{r_i}{r_i+1}\geq r_i \geq \frac{n-i}{i+1}$$
$$\frac{c_i}{c_{i+1}}\frac{{n \choose i+1}}{{n \choose i}}\geq \frac{n-i}{i+1}$$
$$\frac{c_i}{c_{i+1}}\geq 1$$
$$c_i\geq c_{i+1}$$
as required.
\end{proof}

\begin{proposition}
If $\cA$ is as above and $G$ is a graph of order $n$ such that for 
$k=\lceil n/2 \rceil$, 
\begin{gather}
\frac{c_k}{{n \choose k }}\geq \frac{n-k}{k+1}, 
\tag{*}
\end{gather}
then 
the sequence $\{c_i\}$ is unimodal with mode $\lceil n/2 \rceil$.\\
\end{proposition}
\begin{proof}
From Lemma \ref{lemma1}, we have that $c_k\leq c_{k+1}$ for $k< n/2$, and from Lemma \ref{lemma2} we have that $c_k\geq c_{k+1}$ for $k\geq n/2$.
Thus, the sequence $\{c_i\}$ is unimodal with mode $\lceil n/2 \rceil$.\\
\end{proof}

In particular, we have:
\begin{corollary}
\label{cor}
If $\cA$ is as above and $G$ is a graph of order $n$ such that 
for every subset $S\subseteq V(G)$ with $|S|=\lceil n/2 \rceil$, $G[S]\in \cA$,
the sequence $\{c_i\}$ is unimodal with mode $\lceil n/2 \rceil$.
\end{corollary}

\begin{proof}
If there are no subsets $S\subseteq V(G)$ such that $|S|=\lceil n/2 \rceil$ and $G[S]\not\in \cA$, 
then for $k=\lceil n/2 \rceil$, we have $\frac{c_k}{{n \choose k }}=1>\frac{\lfloor n/2 \rfloor}{\lceil n/2 \rceil+1}$.
\\
This can also be shown without using Lemma \ref{lemma2} by noting that in this case we have
$c_i(G) = {n \choose i}$, 
for $i \geq \lceil n/2 \rceil$.
\end{proof}

\subsection{Using random graphs}

We can use Corollary \ref{cor} to show that for many graph properties $\cA$ the  sequence $c_i(G)$
is unimodal for almost all graphs $G$.
In particular, Corollary \ref{cor} applies to all the cases where $\cA$ consists of all 
graphs which contain a fixed induced sugraph $H$.

For the general case we use the following theorem, due to Janson, Luczak, and Ruci\'nski. It was proved in \cite{janson1990exponential}, but a more accessible discussion can be found in \cite{bollobas2011random} (Theorem 4.15) and \cite{frieze2016introduction} (Corollary 23.14).

\begin{theorem}[\cite{janson1990exponential}]
\label{JLR}
For $n\in \NN$, $p\in (0,1)$, let  $G \in \mathcal{G}(n,p)$ be a random graph.
Let $H$ be a fixed graph, and let $X_{n,p}(H)$ be the random variable counting the number of (not induced) subgraphs of $G$ that are isomorphic to $H$. Then:
$$
\log_2 (\PP(X_{n,p}(H)=0)) \leq -C \min \{\EE(X_{n,p}(H'): H' \text{ is a subgraph of H}, e(H')>0\}
$$
for some positive constant $C$.
\end{theorem}
It will be more convenient for us to use the following formulation of Theorem \ref{JLR}, presented as part of the proof of Theorem \ref{JLR} in \cite{janson1990exponential} as equation 3.6:

\begin{proposition}
\label{JLR2}
For $n\in \NN$, $p\in (0,1)$, let  $G \in \mathcal{G}(n,p)$ be a random graph.
Let $H$ be a fixed graph, and let $X_{n,p}(H)$ be the random variable counting the number of (not induced) subgraphs of $G$ that are isomorphic to $H$. Then:
$$
\log_2 (\PP(X_{n,p}(H)=0)) \leq - \left( \sum_{H'} \frac{\#\{H'\subseteq H\}^2}{\EE X_{n,p}(H')} \right)^{-1}
$$
where the sum is over all non isomorphic subgraphs $H'$ of $H$ with at least one edge, and $\#\{H'\subseteq H\}$ is the number of subgraphs of $H$ that are isomorphic to $H'$.
\end{proposition}
We only need the following consequence of the proposition \ref{JLR}:
\begin{corollary}
For $n\in \NN$, $p\in (0,1)$, let  $G \in \mathcal{G}(n,p)$ be a random graph.
Let $H$ be a fixed graph. The probability that $G$ is an $H$-free is bounded by $2^{-cn^{c'}}$ where $c>0$ and $c'\geq 2$
\end{corollary}
\begin{proof}
We evaluate:
$$
\PP(\text{a random graph with $n$ vertices is $H$ free})\leq 
$$
$$
\PP(\text{a random graph with $n$ vertices does not have  $H$ as a subgraph})\leq
$$
\begin{equation}
\label{eq:jrl}
2^{- \left( \sum_{H'} \frac{\#\{H'\subseteq H\}^2}{\EE X_{n,1/2}(H')} \right)^{-1}}
\end{equation}

Using the fact that 
$$\EE X_{n,1/2}(H')= {n \choose |V(H')|}\frac{|V(H')|!}{aut(H')}(1/2)^{|E(H')|}$$ 
where $aut(H')$ is the number of automorphisms of $H'$, we evaluate $  \sum_{H'} \frac{\#\{H'\subseteq H\}^2}{\EE X_{n,1/2}(H')} $: 
$$
\sum_{H'} \frac{\#\{H'\subseteq H\}^2}{\EE X_{n,1/2}(H')}=\sum_{H'}\frac{\#\{H'\subseteq H\}^2\frac{aut(H')2^{|E(H')|}}{|V(H')|!}}{{n \choose |V(H')|}}
$$
Note that the numerator of the fraction in the sum is a positive constant that depends on $H'$, and  asymptotically ${n \choose |V(H')|}\approx n^{|V(H')|}$. Noting that every $H'$ in the sum has at least one edge, and hence at least two vertices, we conclude that there is a $c'\geq 2$ such that
$$
\sum_{H'} \frac{\#\{H'\subseteq H\}^2}{\EE X_{n,1/2}(H')}\leq \frac{c}{n^{c'}}
$$
for some positive constant $c$. Returning to equation \ref{eq:jrl}, we get
$$
\PP(\text{a random graph with $n$ vertices is $H$ free})\leq 2^{-1/c n^{c'}}
$$
as required.
\end{proof}
We only need the following consequence:

\begin{theorem}
\label{thm}
Let $H$ be a fixed graph, and $G\in \mathcal{G}(n,1/2)$ a random graph. Then with high probability, $G$ does not have an $H$ free subgraph with $n/2$ vertices. 
\end{theorem}
\begin{proof}
We bound the probability that a subset of $k$ vertices of $G$ induces an $H$-free graph:
$$\PP(\exists \text{a set of size k in V(G) that induces an H free graph})\leq$$ $$ \EE(\# \text{of sets of size k in V(G) that induces an H free graph})=$$ $${n \choose k}\PP(\text{a random graph with k vertices is H free})\leq $$ $${n \choose k}2^{-ck^{c'}}
\leq \left(\frac{ne}{k}\right)^k 2^{-ck^{c'}}
$$
When $k=n/2$, we have
$$\left(\frac{ne}{k}\right)^k 2^{-ck^{c'}}=(2e)^{n/2}2^{-c(n/2)^{c'}}=(\sqrt{2e})^n 2^{-c(n/2)^{c'}}
$$
Which tends to 0 since $c'\geq 2$.
\end{proof}

Now we are in a position to prove Theorem \ref{th:main-1}: 
\begin{quote}
\em
Let $\cA$ a non-trivial co-hereditary graph property.
Let $c^{\cA}_i(G)$ be the number of subsets $S\subseteq V(G)$ of a graph $G$ with $|S| =i$ such that
$G[S] \in \cA$.
For almost all graphs $G$, the sequence $c^{\cA}_i(G), 0\leq i \leq n(G)$ is unimodal. 
\end{quote}
\begin{proof}
Let $H\in \cA$.
By Theorem \ref{thm} almost all graphs $G$ don't have an induced subgraph with $|V(G)|/2$ vertices that is $H$ free. 
But from Corollary \ref{cor}, for every such graph the sequence $\{c_k\}$ is unimodal with mode $\lceil |V(G)|/2 \rceil$.
\end{proof}

\section{Proof of the Real-rootedness Theorem \ref{th:main-2}}

Now we prove the Real-rootedness Theorem \ref{th:main-2}:

{\bf Theorem \ref{th:main-2} }:
{\em
Let $\cA$ be hereditary and with a graph $G_0 \in \cA$ which is neither a clique nor an edgeless graph.
Then $P_{\cA}(G;x)$ is real-rooted iff $G \in \cA$.
}

Let $\cH$ be a family of graphs  and
$Forb(\cH)$ the class of graphs with no induced subgraph in $\cH$.
We will the following characterization of hereditary graph properties.

\begin{theorem}[Folklore]
A graph property $\cA$ is hereditary iff there is a family of graphs $\cH$ such that
$\cA = Forb(\cH)$.
\end{theorem}

We adapt some definitions from \cite{barton2020acyclic}.

\begin{definition}
Let $\cA$ be a hereditary property, and $G$ a graph.
\begin{enumerate}[(i)]
\item
Denote by $g_{\cA}(G)$ the smallest integer such that there is a set $S\subseteq V(G)$ that does not induce a graph in $\cA$, i.e.
$g_{\cA}(G) = \min\{|S|: S \subseteq V(G), G[S] \not\in \cA \} $
\item
Denote by $g_{\cA}$ the order of the smallest graph that is not in $\cA$. Note that
$g_{\cA} = \min_G \{ g_{\cA}(G)\} $.
\item
Denote by $\triangledown_{\cA}(G)$ the smallest integer such that there is a set $S\subseteq V(G)$ whose complement in $G$ induces a graph in $\cA$, i.e.
$ \triangledown_{\cA}(G) = \min\{|S|: S \subseteq V(G), G[V(G)-S] \in \cA \} $
\end{enumerate}
\end{definition}

Note that for
$ P_{\cA}(G;x) = \sum_{S \subseteq V(G), G[S] \in \cA} x^{|S|}$, 
the  graph polynomial 
$P_{\cA}(G;x)$ is of degree  $n(G) - \triangledown_{\cA}(G)$.

\begin{examples}
\phantom{sdf}
\begin{enumerate}[(i)]
\item
The independence polynomial $I(G;x)$ is of this form with $\cA$ the class of edge-less graphs and $g_{\cA}=2$.
\item
For $\cA$ the class of complete graphs,
$P_{\cA}(G;x) =Cl(G;x)$
is the Clique polynomial, and $g_{\cA}=2$.
\item
Let $\cH$ be a family of graphs with $\mu(\cH)$ the order of its smallest member.
\\
Then $g_{Forb(\cH)} = \mu(\cH)$.
\item
If $\cA$ is hereditary and $g_\cA =2$ then $\cA$ consists either of 
\\
all complete graphs or of all edge-less graphs. 
\end{enumerate}
\end{examples}

\begin{lemma}
\label{lemma:FinalLemma1}
If $G \in \cA$, 
$P_{\cA}(G;x)$ is real rooted.
\end{lemma}

\begin{proof}
In this case
$P_{\cA}(G;x) = \sum_i {n \choose i} x^i = (1+x)^n$.
\\
Therefore $-1 \in \R$ is a root with multiplicity $n$.
\end{proof}
We now consider $P_{\cA}(G;x)$ for $G\not \in \cA$.
\begin{lemma}
\label{lemma:mess}
Let $\cA$ be a hereditary graph property and $G\not \in \cA$ a graph. Then we can write
$$P_{\cA}(G;x) = B(x)  + \left({n \choose g_{\cA}(G)}  - \alpha \right)  x^{g_{\cA}(G)} + C(n)$$
where
$$B = \sum_{i =0}^{g_{\cA}(G)-1} {n \choose i} x^i $$ 
and 
$$C(x) = \sum_{j= g_{\cA}(G)+1}^{n(G) - \triangledown_{\cA}(G)} a_j x^j$$
with $a_j$ an integer for all $j$, and $\alpha$ a positive integer.
\end{lemma}
\begin{proof}
From the definition of $g_{\cA}(G)$, note that for all $0 \leq j < g_{\cA}(G)$, the coefficient of $x^j$ in $P_{\cA}(G;x)$ is ${n \choose i}$, since every set of size $j$ induces a graph in $\cA$. On the other hand, for $j>n(G) - \triangledown_{\cA}(G)$ the coefficient of $x^j$ is $0$, since no set with more then $n(G) - \triangledown_{\cA}(G)$ vertices can induce a graph in $\cA$.
\end{proof}

\begin{lemma}
\label{lemma:FinalLemma2}
Let $f(x) \in \R[x]$ be a polynomial of degree $d$ and $n \geq d$.
\\
$f(x)$ is real-rooted iff $x^n f(\frac{1}{x})$ is real-rooted.
\end{lemma}

\begin{proof}
Let $f(x) = 
x^k \prod_{i=0}^{d-k} (x-c_i)$ with roots $c_i \neq 0$ and $0$ with multiplicity $k$.
\\
Now 
$f(\frac{1}{x}) =
(\frac{1}{x})^k 
\prod_{i=0}^{d-k} (\frac{1}{x}-c_i)$ 
with roots $d_i = \frac{1}{c_i}$.
\\
$x^n f(\frac{1}{x}) = x^{n-k} 
\prod_{i=0}^{d-k} (\frac{1}{x}-c_i)$ 
has root $0$ with multiplicity  $n-k$ 
and the roots $d_i$.
\\
Clearly $d_i$ is real iff $c_i$ is real.
\\
We used that $n \geq d \geq k$.
\end{proof}

\subsection{Using Sturm's Theorem}
Let $F(x) \in \R[x]$ be a real polynomial of degree $n \geq 1$ and leading coefficient $a_n > 0$.
\\
The polynomials
$$F_0(x), F_1(x), F_2(x), \ldots , F_k(x)$$
with respective degrees $d_i$ form
a {\em Sturm sequence} for $F$ if 
\begin{enumerate}[(i)]
\item
$F_0(x) = F(x)$,
\item
$F_1(x) = F(x)'$, the derivative of $F(x)$,
\item
$F_i(x) = -R_{i-2,i-1}(x)$ 
for $i \geq 2$, where
$$
F_{i-2}(x) = F_{i-1}(x) \cdot D_i(x) + R_{i-2,i-1}(x)
$$
where $R_{i-2,i-1}(x)$ is the remainder of the division of $F_{i-2}(x)$ by $F_{i-1}(x)$.
\item
The degrees $d_i$ form a strictly decreasing sequence.
\item
$F_k(x)$ has smallest possible positive degree.
\end{enumerate}
In the paper  \cite{brown2002existence}
the following consequence of Sturm's Theorem is proved
\footnote{The authors say this is stated (imprecisely) in \cite[Page 176]{barbeau2003polynomials}.
They then give a corrected statement and a complete proof.
}:

\begin{theorem}[\cite{brown2002existence}]
\label{thm:FinalThm1}
Let $F(x) \in \R[x]$ be a real polynomial of degree $n \geq 1$ and 
leading coefficient $a_n > 0$.
\\
Let
$F_0(x), F_1(x), F_2(x), \ldots , F_k(x)$ its Sturm sequence with degrees $d_i$.
\\
Then $F(x)$ is real-rooted  iff no $F_i(x)$ has a negative leading coefficient
and $|d_i - d_{i+1}| =1$ for all $j \leq k$.
\end{theorem}


We now compute the degree of the third term in the Sturm sequence for 
$$x^n P_{\cA}(G, \frac{1}{x})$$

\begin{lemma}
\label{le:remainder}
Let $\cA$ be a hereditary property, and let $G \not \in \cA$ be a graph of order $n$.
Let $F(x)=x^n P_{\cA}(G, \frac{1}{x})$, $F_0(x),F_1(x),...,F_k(x)$ its Sturm sequence, and $d_i=\deg(F_i)$.
Then $|d_1-d_2|>1$.
\end{lemma}
\begin{proof}
Let $F_0(x) = x^n P_{\cA}(G, \frac{1}{x})$ and $F_1(x) = F_0'(x)$ its derivative.
To simplify notation we set $g = g_{\cA}(G)$, $\triangledown_{\cA}(G)=\triangledown_{\cA}$ and $d=n-\triangledown_{\cA}$. 
Using the notation of lemma \ref{lemma:mess}, we have 
$$F_0(x)=x^n P_{\cA}(G, \frac{1}{x})=\sum_{j=0}^{g-1} {n \choose j} x^{n-j}+\left({n \choose g}-\alpha \right)x^{n-g}+\sum_{j=g+1}^{d}a_j x^{n-j} $$
For the second term in the Sturm sequence, we have
$$
F_1(x)=F'_0(x)=\sum_{j=0}^{g-1} {n \choose j}(n-j) x^{n-j-1}+\left({n \choose g}-\alpha\right)(n-g)x^{n-g-1}+\sum_{j=g+1}^{d}a_j(n-j) x^{n-j-1}
$$
By applying long division, we can write $F_0(x)=(\frac{1}{n}x+{\frac{1}{n}}) F_1(x)+R(x)$, where 
\begin{gather}
R(x)=\sum_{j=1}^{g-1}\left[{n \choose j}-\frac{n-j}{n}{n \choose j} - \frac{n-j+1}{n}{n \choose {j-1}}\right] x^{n-j}+ 
\label{eq-1} \\
+\left[ \left({n \choose g}-\alpha\right)-\left({n \choose g}-\alpha\right)\frac{n-g}{n}-\frac{n-g+1}{n}{n \choose g-1} \right] x^{n-g}+ \notag \\
+\sum_{j=g+1}^d \left[ a_j-a_j\frac{n-j}{n}-a_{j-1}\frac{n-j+1}{n} \right] x^{n-j}+a_d\frac{n-d}{n}x^{n-d-1} \notag
\end{gather}
The coefficients in Line \ref{eq-1}  vanish, hence
\begin{gather}
R(x) =
\left[\left({n \choose g}-\alpha\right)-\left({n \choose g}-\alpha\right)\frac{n-g}{n}-\frac{n-g+1}{n}{n \choose g-1}\right] x^{n-g}+  \label{eq-2}\\
+\sum_{j=g+1}^d\left[a_j-a_j\frac{n-j}{n}-a_{j-1}\frac{n-k+1}{n}\right] x^{n-k}+a_d\frac{n-d}{n}x^{n-d-1} \notag
\end{gather}
The leading term of $R(x)$ is
$$
\left[\left({n \choose g}-\alpha\right)-\left({n \choose g}-\alpha\right)\frac{n-g}{n}-\frac{n-g+1}{n}{n \choose g-1}\right]x^{n-g} = \alpha(\frac{n-g}{n}-1)x^{n-g}
$$
which has degree $d_2=n-g < n-2$, and since $d_1= n-1$ the claim follows.
\end{proof}
Using the above lemmas we can finally prove Theorem \ref{th:main-2}:

{\bf Theorem \ref{th:main-2} }:
{\em
Let $\cA$ be hereditary and with a graph $G_0 \in \cA$ which is neither a clique nor an edgeless graph.
Then $P_{\cA}(G;x)$ is real-rooted iff $G \in \cA$.
}
\begin{proof}
If $G\in \cA$, then $P_{\cA}(G;x)$ is real-rooted by Lemma \ref{lemma:FinalLemma1}. On the other hand, if $G\not \in \cA$, from Lemma \ref{lemma:FinalLemma2} we have that $P_{\cA}(G;x)$ is real-rooted if and only if $x^nP_{\cA}(G;1/x)$ is real rooted, and from Lemma \ref{le:remainder} and Theorem \ref{thm:FinalThm1}, we have that $x^nP_{\cA}(G;1/x)$ is not real rooted, so we are done.
\end{proof}
\section{Conclusions and open problems}
In Theorem \ref{th:main-1} we have shown that the
the generating function $P_{\cA}(G;x)$ of a co-hereditary graph property $\cA$ 
is unimodal for almost all graphs.

\begin{question}
Under what conditions can almost unimodality be improved to unimodality?
\end{question}

\begin{question}
Under what conditions can unimodality be improved to log-concavity?
\end{question}

It also follows from Theorem \ref{th:main-1} that for $\cA$ hereditary the graph parameter
$$
{n \choose i} - c_i
$$
is unimodal.  The sequence $n \choose i$ is known to be log-concave.
However, the difference between a log-concave sequence and a unimodal sequence
need not be unimodal,
even if they have the same length and mode.

\begin{question}
What more can we say about $c_i$ if $\cA$ is hereditary?
\end{question}

In particular:

\begin{question}
Is the independence polynomial unimodal for almost all graphs?
\end{question}

With Theorem \ref{th:main-2} we characterized real-rootedness of
$P_{\cA}(G;x)$ for hereditary graph properties $\cA$.
This leaves several questions open:

\begin{problem}
Characterize the hereditary graph properties $\cA$ such that
$P_{\cA}(G;x)$ is unimodal, although not necessarily real-rooted.
\end{problem}

An augmented graph property $\cP$ is {\em downward monotone} if for every $(G, S) \in \cP$,
for every induced subgraph of $G'<G$ and subset $S'\subseteq S$ such that $S'\subseteq V(G')$, $(G',S')\in \cP$

\begin{problem}
Can Theorem \ref{th:main-2} be extended to (downward monotone) augmented graph properties $\cQ$:
$F_{\cA}(G;x)$ is real-rooted iff $G \in \cA_{\cQ}$ for some suitable graph property $\cA_{\cQ}$?
\end{problem}

\begin{problem}
Characterize the hereditary graph properties $\cA$ such that
$F_{\cQ}(G;x)$ is unimodal, although not necessarily real-rooted.
\end{problem}

\newcommand{\etalchar}[1]{$^{#1}$}

\end{document}